\documentclass[final]{siamltex}         

\usepackage[dvips]{color}

\usepackage{graphicx}

\usepackage{subfigure}

\graphicspath{{C:/NS/Pictures/}}

\title{An effective approach to the solution of a system of nonlinear differential equations in partial derivatives.}

\author{ A. Tsionskiy, M. Tsionskiy  }               

\begin{document}             

\maketitle                   

\begin{abstract}
There are few approaches to the solution of a system of nonlinear differential equations in partial derivatives, for example $\cite{NK87} - \cite{EK98}$. In our paper we propose an approach that was used to solve  the Navier-Stokes equations in three dimensional space. This solution is described in details in article "Existence, uniqueness and smoothness of solution for 3D Navier-Stokes equations with any smooth initial velocity" $\cite{TT12}$. The authors expect that it can be successfully applied to other systems of nonlinear differential equations in partial derivatives.

\end{abstract}



\pagestyle{myheadings}
\thispagestyle{plain}
\markboth{A. TSIONSKIY, M. TSIONSKIY}{AN EFFECTIVE APPROACH TO THE SOLUTION}

\section{The mathematical setup}\ 
$\\$
$\\$
The system of non-linear partial differential equations is given by

\begin{equation}\label{eqn1}
F_{k, l}(x, \vec{u}, D^{\alpha}\vec{u})\;+\;F_{k, n}(x, \vec{u}, D^{\beta}\vec{u})\; =\;\;f_{k}(x)\;\;\;\;\;(x \in R^{N},\;\;{1\leq k \leq N}) 
\end{equation}

Here \(\vec{u}(x)=(u_{k}(x)) \in R^{N},\;\; ({1\leq k \leq N})  \;-\; \)is an unknown vector-function;  \(\; f_{k}(x)\;\)are components of a given, externally applied force \(\vec{f}(x)\); \(F_{k, l}(x, \vec{u}, D^{\alpha}\vec{u})\) are linear parts of equations; \(F_{k, n}(x, \vec{u}, D^{\beta}\vec{u})\) are non-linear parts of equations; \(\alpha = (\alpha_{1}, \ldots, \alpha_{n})\) is a multi-index of non-negative integers \(\alpha_{1}, \ldots, \alpha_{n}\), and \(\beta = (\beta_{1}, \ldots, \beta_{n})\) is a multi-index of non-negative integers \(\beta_{1}, \ldots, \beta_{n}\). Then we have \(D^{\alpha} = D_{1}^{\alpha_{1}} \ldots D_{n}^{\alpha_{n}},\; D^{\beta} = D_{1}^{\beta_{1}} \ldots D_{n}^{\beta_{n}}\), where \(D_{i} = \frac{\partial }{\partial x_{i}}\;\;(1\leq i \leq N).\) The order of (\ref{eqn1}) is defined as the highest order of a derivative occurring in the system of equations. 

System of equations (\ref{eqn1}) can be considered in the whole space  or in some sub domain of it. In the first case the definition of the solution space has to include conditions on the behavior of the solutions at infinity. In case if system (\ref{eqn1}) is considered in some sub domain of the whole space one or more boundary conditions have to be imposed. For this approach the boundary conditions (or conditions at infinity) must include only linear operators. 

A system of non-linear partial differential equations together with boundary conditions gives a non-linear problem, which must be considered in an appropriate function space. 

\section{Integration of linear part of the system of partial differential equations (\ref{eqn1})}\ 
$\\$
$\\$
To start the process of solution of a non-linear problem let us add (\( - F_{k, n}(x, \vec{u}, D^{\beta}\vec{u})\;\)) to both sides of the equations (\ref{eqn1}). Then we have:

\begin{equation}\label{eqn7}
F_{k, l}(x, \vec{u}, D^{\alpha}\vec{u})\; =\;\;f_{k}(x) -\;F_{k, n}(x, \vec{u}, D^{\beta}\vec{u})\;\;\;\;\;\;(x \in R^{N},\;\;{1\leq k \leq N}) 
\end{equation}

Let us denote

\begin{equation}\label{eqn13}
\tilde{f}_{k}(x,t)\; = \; f_{k}(x,t) \; - \; \;F_{k, n}(x, \vec{u}, D^{\beta}\vec{u})\;\;\;\;\;\;({1\leq k \leq N})
\end{equation}

or we can present it in the vector form:

\begin{equation}\label{eqn14}
\vec{\tilde{f}}(x,t)\; = \; \vec{f}(x,t) \; - \;\vec{F_{n}}(x, \vec{u}, D^{\beta}\vec{u})
\end{equation}
\\
Then we should solve the system of linear partial differential equations. Corresponding methods of this integration are presented in V.P. Palamodov $\cite{VP70}$, G.E. Shilov $\cite{gS01}$, L. Hormander $\cite{LH83}$, S. Mizohata $\cite{SM73}$, J.F. Treves $\cite{JFT61}$, R. Courant, D. Hilbert  $\cite{RCD89}$.

Result of this stage is the system of integral equations for vector-function $\vec{u}$. 
Non-linear parts  $(\ref{eqn14})$ are integrands for integrals in  the system of integral equations.

We have to show equivalence of the problem in differential form $(\ref{eqn1})$ and in the form of an appropriate system of integral equations.

\section{The fixed point principle}$\cite{KA64},\;$ $\cite{VT80},\;$ $\cite{WR73},\;$$\; \cite{KS01},\;$$\; \cite{GD03},\;$$\;\cite{ADL97}\;$
$\\$
$\\$
We can use the fixed point principle (L.V. Kantorovich, G.P. Akilov $\cite{KA64},\;$ V.A. Trenogin  $\cite{VT80},\;$W. Rudin $\cite{WR73},\;$W.A. Kirk and B. Sims$\; \cite{KS01},\;$A. Granas and J. Dugundji$\; \cite{GD03},\;$J.M. Ayerbe Toledano, T. Dominguez Benavides, G. Lopez Acedo$\;\cite{ADL97}\;$) for the solution of received equivalent system of integral equations.
$\\$
$\\$
This approach was used in the article "Existence, uniqueness and smoothness of solution for 3D Navier-Stokes equations with any smooth initial velocity" $\cite{TT12}$
$\\$
$\\$


\begin{thebibliography}{2}

\bibitem {NK87}
  {\sc N.V. Krylov},
  \emph{Nonlinear elliptic and parabolic equations of the second order,}
 Reidel, 1987, (Translated from Russian).
\bibitem {PR64}
  {\sc P.R. Garabedian},
  \emph{Partial differential equations,}
 Wiley, 1964.
\bibitem {JS07}
  {\sc D. W. Jordan, P. Smith},
  \emph{Nonlinear Ordinary Differential Equations ,}
 (fourth ed.). Oxford Univeresity Press, 2007.
\bibitem {KH01}
  {\sc Khalil, Hassan K.},
  \emph{Nonlinear Systems,}
 Prentice Hall, 2001. 
\bibitem {EK98}
  {\sc Erwin Kreyszig},
  \emph{Advanced Engineering Mathematics,}
 Wiley, 1998.    
\bibitem {oL69}
  {\sc O. Ladyzhenskaya},
  \emph{The Mathematical Theory of Viscous Incompressible Flows,}
 (2nd edition), Gordon and Breach, 1969.
\bibitem {jL34}{\sc J. Leray},
\emph { Sur le Mouvement d'un  Liquide Visquex Emplissent l'Espace.}
 Acta Math. J. 63 (1934), 193-248.
\bibitem {TT10}{\sc A. Tsionskiy, M. Tsionskiy},
\emph { Solution of the Cauchy problem for the Navier - Stokes and Euler equations,}
  arXiv:1009.2198v3, 2010.
\bibitem {TT11}{\sc A. Tsionskiy, M. Tsionskiy},
\emph { Research of convergence of the iterative method for solution of the Cauchy problem for the Navier - Stokes equations based on estimated formula ,}
  arXiv:1101.1708v2, 2011.
\bibitem {TT12}{\sc A. Tsionskiy, M. Tsionskiy},
\emph { Existence, uniqueness and smoothness of solution for 3D Navier-Stokes equations with any smooth initial velocity,}
   arXiv:1201.1609v6, 2012.     
\bibitem {VP70}{\sc V. P. Palamodov},
\emph { Linear Differential Operators with Constant Coefficients.}
 Berlin, New York, Springer-Verlag, 1970.      
\bibitem {gS01}{\sc G.E. Shilov},
  \emph { Elementary functional analysis,}
  Cambridge, Mass.: MIT Press, 1974. 
\bibitem {LH83}{\sc L. Hormander},
\emph { The Analysis of Linear Partial Differential Operators I - IV.}
 Berlin, New York, Springer Verlag, 1983 - 1985. 
\bibitem {SM73}{\sc S. Mizohata},
\emph { The Theory of Partial Differential Equations.}
 Cambridge Univ. Press, 1973. 
\bibitem {JFT61}{\sc J.F. Treves},
\emph { Lectures on linear partial differential equations with constant coefficients.}
 Rio de Janeiro : Instituto de Matema´tica Pura e Aplicada do Conselho Nacional de Pesquisas, 1961. 
\bibitem {RCD89}{\sc R. Courant, D. Hilbert},
\emph { Methods of Mathematical Physics.}
Volume 2, Wiley-VCH; 1989.   
\bibitem {GC68}{\sc I.M. Gel'fand, G.E. Chilov},
\emph { Generalized functions./Volume 2, Spaces of fundamental and generalized functions.}
 New York; London: Academic Press, 1968. 
\bibitem {KA64}{\sc L.V. Kantorovich, G.P. Akilov},
\emph { Functional analysis in normed spaces.}
  Oxford/London/Edinburgh/New York/Paris/Frankfurt, Pergamon Press., 1964.  
\bibitem {VT80}{\sc V.A. Trenogin},
\emph { Functional'nyiy analiz (Functional analysis).}
  Nauka, Moskva, GRFML, Russian, 1980. 
\bibitem {WR73}{\sc W. Rudin},
\emph { Functional analysis.}
  New York St. Louis San Francisco Dusseldorf Johannesburg Kuala Lumpur London Mexico Montreal New Delhi Panama Rio de Janeiro Singapore Sydney Toronto, McGRAW-HILL BOOK COMPANY, 1973.  
\bibitem {KS01}{\sc W.A. Kirk and B. Sims},
\emph { Handbook of Metric Fixed Point Theory,}
  Kluwer Academic, London, 2001. 
\bibitem {GD03}{\sc A. Granas and J. Dugundji},
\emph { Fixed Point Theory,}
  Springer-Verlag, New York, 2003. 
\bibitem {ADL97}{\sc J.M. Ayerbe Toledano, T. Dominguez Benavides, G. Lopez Acedo},
\emph { Measures of Noncompactness in Metric Fixed Point Theory,}
  Birkhauser Verlag, Basel - Boston - Berlin, 1997.
\bibitem {RR64}{\sc A.P. Robertson, W.J. Robertson},
\emph { Topological Vector Spaces.}
 Cambrige University Press, 1964.
\bibitem {RN72}{\sc F. Riesz, B.S - Nagy},
\emph { Lecons D'ANALYSE FONCTIONNELE.}
 Academiai Kiado, Budapest 1972, 6 edition.
 
{\sc  F. Riesz, B.Szokefalvi - Nagy},
\emph { Functional analysis.}
 New York: Ungar., 1972.     
\bibitem {RM76}{\sc R.H. Martin},
\emph { Nonlinear Operators and Differential Equations in Banach Spaces,}
  John Wiley and Sons, 1976.         
\end{thebibliography}
\end{document}